\documentclass[10pt]{amsart}

\usepackage{amssymb}
\usepackage{amsmath}

\newtheorem{theorem}{Theorem}[section]
\newtheorem{lemma}[theorem]{Lemma}

\newtheorem{corollary}[theorem]{Corollary}

\newcommand{\p}[1]{\textrm{#1}}

\theoremstyle{definition}

\theoremstyle{remark}
\newtheorem{remark}[theorem]{Remark}

\numberwithin{equation}{section}

\begin{document}

\title{A deterministic version of Pollard's $p-1$ algorithm}

\author{Bartosz \'Zra\l ek}

\address{Institute of Mathematics, Polish Academy of Sciences, 00-956 Warsaw, Poland}

\email{b.zralek@impan.gov.pl}
\subjclass[2000]{Primary 11Y16; Secondary 11Y05, 68Q10}


\keywords{Pollard's $p-1$ method, derandomization, Euler's $\varphi$-function and factorization}

\begin{abstract}
In this article we present applications of smooth numbers to the unconditional derandomization 
of some well-known integer factoring algorithms.

We begin with Pollard's $p-1$ algorithm, which finds in random polynomial time the prime divisors 
$p$ of an integer $n$ such that $p-1$ is smooth. We show that these prime factors can be recovered in deterministic polynomial time. We further generalize this result to give a partial derandomization of the 
$k$-th cyclotomic method of factoring ($k\ge 2$) devised by Bach and Shallit. 

We also investigate reductions of factoring to computing Euler's totient function $\varphi$.
We point out some explicit sets of integers $n$ that are completely factorable in deterministic polynomial time given $\varphi(n)$. These sets consist, roughly speaking, of products of primes $p$ satisfying, with the exception of at most two, certain conditions somewhat weaker than the smoothness of $p-1$.
Finally, we prove that $O(\ln n)$ oracle queries for values of $\varphi$ are sufficient to completely
factor any integer $n$ in less than $\exp\Bigl((1+o(1))(\ln n)^{\frac{1}{3}} (\ln\ln n)^{\frac{2}{3}}\Bigr)$
deterministic time.
\end{abstract}

\maketitle

\section{Introduction}

A fundamental question of algorithmic number theory, in particular, and complexity theory, in general, asks whether there are computational problems which cannot be solved efficiently without the use of randomness. If the answer is no, then we would say that every algorithm can be derandomized. The issue surely has a philosophical flavour, but above all is essential for the development of mathematics. As a rule, derandomization presupposes  making the most of the rich mathematical structures involved. 
It gives rise to new ideas, subtle refinements of existing ones, or, in the worst case, generates fascinating open problems. One of these problems, determining the complexity of primality testing, has been brilliantly solved in \cite{aks}:  primes are recognizable in deterministic polynomial time.

In this article we present applications of smooth numbers to the unconditional derandomization of some well-known integer factoring algorithms. Recall that a smooth number is a product of small primes 
(small relative to, say, $n$ meaning polynomial in the size of $n$).

In sections 3 and 4 we analyze Pollard's $p-1$ method \cite{pollard}, important both in theory and practice \cite{rsa,elliptic}. Pollard's algorithm finds in random polynomial time those prime divisors $p$ 
of an integer $n$ for which $p-1$ is smooth. We show that such prime factors can be recovered in deterministic polynomial time (corollary \ref{pollard}). Let us merely indicate the two ingredients of the proof. The first comes from F\"{u}rer \cite{furer}, Fellows and Koblitz \cite{fellows}, and also Konyagin and Pomerance \cite{konyagin}: take small integers or, what amounts to the same, small primes to generate a large subgroup $G$ of $\mathbb{Z}_n^*$. The second is a novel idea inspired by the Pohlig-Hellman algorithm \cite{pohlig} for computing discrete logarithms. Namely, let $H$ be the group generated by two elements $a$ and $b$ of $\mathbb{Z}_n^*$, both having smooth order. Then given $a$, $b$ and their orders, we can compute a generator of $H$ or a nontrivial divisor of $n$ in deterministic polynomial time. This result is easily extended by induction to any number of given generators for $H$ (corollary \ref{cykliczna}). We apply it with $H=G$.

In section 5 we give a partial derandomization of the $k$-th cyclotomic method of factoring devised by Bach and Shallit \cite{cyclo}. This method is used to find in random polynomial time such prime factors $p$ of an integer $n$ that the value at $p$ of the $k$-th cyclotomic polynomial is smooth. For the reader's convenience, we first treat the simpler case $k=2$ (theorem \ref{williams}), corresponding to Williams' $p+1$ method \cite{williams}, then that of an arbitrary $k$, $k\ge 2$ (theorem \ref{cyclotomic}). The arguments involve more than the derandomization of the $p-1$ algorithm: some elementary algebraic number theory and a lemma proved in \cite{zralek}.

In the last three sections, we attempt to make some progress on a famous open problem: is factoring reducible in deterministic polynomial time to computing Euler's totient function $\varphi$? (cf. problem 23 of \cite{adleman})

In section 6 we discuss the current state of the art. Miller \cite{miller} found a reduction whose correctness depends on the Extended Riemann Hypothesis (ERH). Rabin \cite{rabin} obtained an unconditional reduction at the cost of giving up determinism. A relatively recent result of Burthe \cite{burthe} yields a reduction for almost all integers, but these cannot be simply described. 

In section 7 we point out some explicit sets of integers $n$ that are completely factorable in deterministic polynomial time given $\varphi(n)$ (theorem \ref{graf phi}). These sets consist, roughly speaking, of products of primes $p$ satisfying, with the exception of at most two, certain conditions somewhat weaker than the smoothness of $p-1$.   

In section 8 we study the deterministic complexity of factoring given an oracle for the function $\varphi$. Suppose that we want to factor into primes the integer $n$. Our idea is first to query the oracle for the iterations $\varphi(n)$, $\varphi^2(n)$, $\varphi^3(n)$, etc. until $\varphi^k(n)=1$. Then to come back up
to the complete factorization of $n$ ($n=\varphi^0(n)$) by a recursive procedure, which recovers the prime factorization of $\varphi^{l-1}(n)$ from the prime factorization of $\varphi^l(n)$, starting with $l=k$. We are basically left with the task of finding the prime factorization of an integer $n$ given the complete factorization of $\varphi(n)$. In the hard case, all the prime divisors of $n$ are congruent to $1$ modulo a large integer $A$ that we compute; we further retrieve the missing information either by a direct search or by factoring the polynomial whose coefficients are the coefficients of $n$ in base $A$ 
(lemma \ref{poly}). The resulting algorithm runs in less than 
$\exp \Bigl((1+o(1))(\ln n)^{\frac{1}{3}} (\ln\ln n)^{\frac{2}{3}}\Bigr)$ deterministic time (theorem \ref{faktoryzacja}). Consequently, factoring is reducible in deterministic subexponential time to computing $\varphi$ (corollary \ref{podwykladniczy}).

\section{Notation}

\noindent
Throughout the text $n$ is an odd integer, and $p,q,s$ are prime numbers.
\newline
The greatest common divisor, respectively the least common multiple, of the integers $a,b$
is denoted by $(a,b)$, respectively $\p{LCM}(a,b)$.
\newline
We let $v_s(m)$ be the exponent of the highest power of $s$ dividing $m$.
\newline
For $G$ a group, $\mathcal{B}\subset G$, $b\in G$, we should denote by 
$\langle\mathcal{B}\rangle_G$ the subgroup of $G$ generated by $\mathcal{B}$, 
and denote by $\p{ord}_G(b)$ the order of $b$ in $G$. However, if $G=\mathbb{Z}_d^*$,
respectively $G=\mathbb{Z}_d[\sqrt{m}]^*$, we will just write $\langle\mathcal{B}\rangle_d$
and $\p{ord}_d(b)$, respectively  $\langle\mathcal{B}\rangle_{d,m}$ and $\p{ord}_{d,m}(b)$.
\newline
The cyclic group with $m$ elements is denoted by $C_m$.
\newline
The symbol $\mathbb{P}$ stands for the set of all prime numbers.
We denote by $p_{-}(m)$, respectively $p_{+}(m)$, the least, respectively the largest, prime 
dividing $m$.
\newline
We use $a_i$ to represent the $i$-th coordinate of 
$a\in\mathbb{Z}_n^*=\bigoplus\limits_{q\mid n}\mathbb{Z}_{q^{v_q(n)}}^*$.
\newline
We recall the definitions of the familiar number-theoretic functions appearing in the text:
\begin{itemize}
\item[] $\varphi(m)=\#\lbrace d\le m:~(d,m)=1\rbrace$ (Euler's totient function),
\item[] $\omega(m)=\sum\limits_{p\mid m}1$ and $\Omega(m)=\sum\limits_{p\mid m}v_p(m)$,
\item[] $\psi(x,y)=\#\lbrace m\le x:~p_{+}(m)\le y\rbrace$.
\end{itemize}
We will make frequent use of the following theorem proved in \cite{konyagin}:

\begin{theorem}[Konyagin, Pomerance]                   
\label{konyaginpomerance}
If $n\ge 4$ and $2\le(\ln n)^c\le n$, then $\psi(n,(\ln n)^c)>n^{1-\frac{1}{c}}$.
\end{theorem}

\noindent
We will always assume that its hypotheses are satisfied when $c$ is fixed 
(this is natural in the task of factoring $n$).
In the last section another estimation of $\psi$ will be applied.

\begin{theorem}[Canfield et al.]
\label{canfield}
There is an effective, positive constant $C$ such that for $x, y\ge 1$
and $u:=\frac{\ln x}{\ln y}\ge 3$ we have
\begin{equation*}
\psi(x, y)\ge x\exp
\left[-u\left\{\ln(u\ln u)-1+\frac{\ln\ln u-1}{\ln u}+C\left(\frac{\ln\ln u}{\ln u}\right)^2\right\}\right].
\end{equation*}
\end{theorem}

\section{Pollard's $p-1$ factoring algorithm}
\label{standardowypollard}

We first sketch the ideas behind the probabilistic version of Pollard's $p-1$ factorization method.
Let $n$ be an odd integer, not a prime power.
Assume that we are given an integer $M$ such that $p-1\mid M$ for some $p\mid n$ (for the moment 
we do not consider the issue of finding a suitable $M$). Choose $b\in\mathbb{Z}_n^*$. 
By Fermat's little theorem we have $b^M\equiv 1 (p)$ and thus $d:=(b^M-1,n)>1$. If additionally $d<n$, 
then $d$ is a nontrivial divisor of $n$. But what if $d=n$, i.e. $b^M=1$? We can pick another 
element of $\mathbb{Z}_n^*$. We can also hope to find a nontrivial factor of $n$ in the sequence 
$(b^{\frac{M}{2^l}}-1,n)_{l=1,\ldots,v_2(M)}$, as all square roots of $1$ in $\mathbb{Z}_n^*$ are 
of the form 
$(\pm 1,\ldots,\pm 1)\in\mathbb{Z}_n^*=\bigoplus\limits_{q\mid n}\mathbb{Z}_{q^{v_q(n)}}^*$.
It turns out that the expected number of random $b\in\mathbb{Z}_n^*$ needed to split $n$ does not exceed $2$.

\begin{theorem}[Rabin]
\label{probabilistycznypollard}
Let $n$ be odd, $n>2$, $M$ be even,
\begin{equation*}
\mathcal{F}(M)=\{b\in\mathbb{Z}_n^*:~b^M\ne 1\},
\end{equation*}
\begin{equation*}
\mathcal{S}(M)=\{b\in\mathbb{Z}_n^*\setminus\mathcal{F}(M):~
\exists_{1\le l\le v_2(M)}~1<(b^{\frac{M}{2^l}}-1,n)<n\}.
\end{equation*}
Then {{\Large $\frac{\#(\mathcal{F}(M)\cup\mathcal{S}(M))}{\varphi(n)}$}~~$\ge~~1-2^{1-\omega(n)}$}.
\end{theorem}

Note that we want not only $M$ to be a multiple of $p-1$ for some (a priori unknown) $p\mid n$, 
but also $\ln M$ to be relatively small (e.g., bounded by a fixed power of $\ln n$), so that raising to the 
power $M$ (or $\frac{M}{2^l}$) modulo $n$ does not take too much time. Suppose that $n$ has 
a prime divisor $p$ such that $p-1$ is smooth, say $p_{+}(p-1)\le(\ln n)^u$. Set
$M=\prod\limits_{q\le(\ln n)^u}q^{\bigl[\frac{\ln n}{\ln q}\bigr]}$.
Then $M$ satisfies the two conditions, since
$\ln M\le\sum\limits_{q\le(\ln n)^u}\frac{\ln n}{\ln q}\ln q=\pi((\ln n)^u)\ln n=
O${\Large$(\frac{(\ln n)^{u+1}}{u\ln\ln n})$}
from Chebyshev's theorem. By contrast, there is no efficient method of finding $M$ if $n$ is
not divisible by a prime $p$ as above.
\newline
As before suppose that $n$ is odd, divisible by at least two different primes $p$ and $q$.
It is well known that if a multiple $M$ of $p-1$ is given, then the previously described
search for a nontrivial factor of $n$ can be derandomized under the ERH. Without loss of
generality assume that $b^M\equiv1 (n)$ for all $b<2(\ln n)^2$. 

\begin{theorem}[Bach]
\label{ankeny}
Suppose that the ERH is true.
Let $n\ge 3$, $\chi$ be a nonprincipal character modulo $n$. There is an integer $b<2(\ln n)^2$
such that $\chi(b)\ne 1$.    
\end{theorem}

Using this theorem, we can easily prove the existence of $b<2(\ln n)^2$ such that 
for some $l$, $b^{\frac{M}{2^l}}-1$ is divisible by $q$ or $p$, but not both. 
We apply it  with $\chi$ induced by the quadratic character
$\left( \frac \cdot p \right), \left( \frac \cdot q \right), \left( \frac \cdot {pq} \right)$ 
when $v_2(p-1)>v_2(q-1), v_2(p-1)<v_2(q-1), v_2(p-1)=v_2(q-1)$, respectively.

\section{A deterministic variant of Pollard's $p-1$ factoring algorithm}
\label{jeden}

Our basic framework is as follows. Let ${\mathcal{B}}=\{2,3,\ldots,[(\ln n)^2]\}$. Assume that
we are given an integer $M$ together with its complete factorization such that
$b^M\equiv 1 (n)$ for every $b\in\mathcal{B}$. We want to find a simple and not restrictive condition
on $n$ under which $n$ is factorable in deterministic polynomial time in $\ln n$ and $\ln M$.
The starting point is a reformulation of the primality criterion from \cite{fellows}.
We restate the argument for completeness  and clarity of exposition.

\begin{theorem}[Fellows-Koblitz]
\label{F-K}
Let ${\mathcal{B}}=\{2,3,\ldots,[(\ln n)^2]\}$, $\mathcal{B}\subset\mathbb{Z}_n^*$.
Then $n$ is prime if and only if the following conditions are satisfied.
\renewcommand{\labelenumi}{(\roman{enumi})}
\begin{enumerate}
\item
$\p{ord}_p(b)=\p{ord}_n(b)$ for every $b\in\mathcal{B}$ and $p\mid n$.
\item
$\p{LCM}_{b\in\mathcal{B}}(\p{ord}_n(b))>\sqrt{n}$. 
\end{enumerate}
\begin{proof}
Suppose $n$ is prime. Condition (i) is then a tautology. We check condition (ii).
The group $\langle\mathcal{B}\rangle_n$ is cyclic, since $n$ is prime. Therefore
\begin{equation*}
\p{LCM}_{b\in\mathcal{B}}(\p{ord}_n(b))=\#\langle\mathcal{B}\rangle_n\ge\psi(n,(\ln n)^2)>\sqrt{n},
\end{equation*}
where the last inequality follows from theorem \ref{konyaginpomerance}.
\newline
Assume now that conditions (i) and (ii) are satisfied. Let $p=p_{-}(n)$. We then have 
$\p{ord}_p(b)=\p{ord}_n(b)$ for all $b\in\mathcal{B}$ and thus
\begin{equation*}
\p{LCM}_{b\in\mathcal{B}}(\p{ord}_p(b))=\p{LCM}_{b\in\mathcal{B}}(\p{ord}_n(b))>\sqrt{n}.
\end{equation*}
However $\p{LCM}_{b\in\mathcal{B}}(\p{ord}_p(b))\mid p-1$. Consequently
$p>\sqrt{n}$; hence $n\in\mathbb{P}$.
\end{proof}
\end{theorem}
Let $b\in\mathbb{Z}_n^*$, $p\mid n$. Recall that $\p{ord}_p(b)<\p{ord}_n(b)$ is equivalent to
$p\mid b^{\frac{\p{ord}_n(b)}{s}}-1$ for some $s\mid\p{ord}_n(b)$. 
If $(b^{\frac{\p{ord}_n(b)}{s}}-1,n)>1$ for some $s\mid\p{ord}_n(b)$, then we will say that $b$ is
a \textit{Fermat-Euclid witness for} $n$. 
Checking conditions (i) and (ii) therefore reduces to factoring the orders of the elements of 
$\mathcal{B}$, which can be done efficiently under our assumption on $M$.
Taking $M=n-1$ yields a deterministic polynomial time algorithm for deciding the primality of integers
$n$ such that $n-1$ is smooth. Actually, a stronger test, in which only a part of $n-1$ exceeding
$n^{\frac{1}{2}+\varepsilon}$ ($\varepsilon>0$) is assumed to be smooth, was first discovered by F\"{u}rer \cite{furer}.  Konyagin and Pomerance \cite{konyagin} further reduced the exponent 
$\frac{1}{2}+\varepsilon$ to $\varepsilon$. The key point is that beside searching some other appropriately chosen ``small'' subset $\mathcal{B}$ of $\mathbb{Z}_n^*$ for Fermat-Euclid witnesses for $n$, one can also check the cyclicity of $\langle\mathcal{B}\rangle_n$. The authors verify this stringent condition by applying the classic Pohlig-Hellman technique \cite{pohlig} of discrete logarithm computation in a prime field. Here we will in a sense extend this technique for the purpose of splitting the integer $n$.
\newline
Suppose for greater generality that $\mathcal{B}$ is any subset of $\mathbb{Z}_n^*$.
We will describe below a deterministic algorithm that finds a generator of 
$\langle\mathcal{B}\rangle_n$ or, particularly in the case when $\langle\mathcal{B}\rangle_n$
is not cyclic, a nontrivial divisor of $n$. This algorithm runs in polynomial time if $\mathcal{B}$
consists of elements having smooth orders in $\mathbb{Z}_n^*$. By induction, it is sufficient to restrict our attention to the case $\#\mathcal{B}=2$, say $\mathcal{B}=\lbrace a,b\rbrace$. 
\newline
We assume temporarily that $\p{ord}_n(a)=s^v$, $b^{s^v}=1$ with $s\in\mathbb{P}$,
$v\in\mathbb{N}$. Let $n=p_1^{e_1}\cdot\ldots\cdot p_k^{e_k}$ be the complete factorization of $n$.
There exist an $i$, $1\le i\le k$, such that $\p{ord}_{p_i^{e_i}}(a_i)=s^v$.
Since $b_i^{s^v}=1$ and $\mathbb{Z}_{p_i^{e_i}}^*$ is cyclic, we have $a_i^l=b_i$ for some uniquely determined, less than $s^v$, natural number $l$. Write $l$ in base $s$: $l=\sum\limits_{0\le r<v}l_rs^r$.
Set $l_{-1}=0$ and reason by induction. Assume we have computed $l_{-1},\ldots,l_m$, where
$-1\le m\le v-2$. Put $c=ba^{-\sum\limits_{-1\le r\le m}l_rs^r}$. Then $c_i=a_i^{\sum\limits_{m<r<v}l_rs^r}$.
Therefore $c_i^{s^{v-m-2}}=a_i^{l_{m+1}s^{v-1}}$. Denote $(c^{s^{v-m-2}}-a^{js^{v-1}},n)$ by $d_j$.
We successively compute $d_0,d_1,\ldots$, until we get $d_j>1$ for some $j\le s-1$. This will happen,
because $p_i^{e_i}\mid d_{l_{m+1}}$. If moreover $d_j<n$, then $d_j$ is a nontrivial factor of $n$.
Otherwise, $d_j=n$. In particular, $c_i^{s^{v-m-2}}=a_i^{js^{v-1}}$. Hence $j=l_{m+1}$. Eventually, if
$m=v-2$, then $d_{l_{m+1}}=n$ implies $b=a^l$. More formally we use the ensuing algorithm.
\newline

\noindent
PH($n,a,b,s,v,w$)
\{$a,b\in\mathbb{Z}_n^*,s\in\mathbb{P},\p{ord}_n(a)=s^v,\p{ord}_n(b)=s^w$\}
\begin{enumerate}
\item
If $w>v$ then interchange $a$ and $b$
\item
For $j=1$ to $s-1$ compute $a^{js^{v-1}}$
\item
Let $c=b$
\item
For $m=-1$ to $v-2$ do
\begin{enumerate}
\item
Let $j=0$
\item
While $(c^{s^{v-m-2}}-a^{js^{v-1}},n)=1$ do $j=j+1$
\item
Let $d=(c^{s^{v-m-2}}-a^{js^{v-1}},n)$. If $d\ne n$ then return $d$
\item
Let $c=ca^{-js^{m+1}}$
\end{enumerate}
\end{enumerate}

\begin{theorem}
\label{PH}
Let $a,b\in\mathbb{Z}_n^*,s\in\mathbb{P},\p{ord}_n(a)=s^v,\p{ord}_n(b)=s^w$. If the algorithm
PH($n,a,b,s,v,w$) does not find a nontrivial divisor of $n$, then $\langle a,b\rangle_n$ is cyclic.
This algorithm uses $O((s+u\ln s)u(\ln n)^2)$ operations, where $u=\max(v,w)$.
\begin{proof}
The correctness of PH($n,a,b,s,v,w$) follows from the preceding discussion. 
Step 2 requires $O(u(\ln n)^2\ln s+s(\ln n)^2)$ operations.
The total number of operations used by step 4b in the loop 4 is $O(u^2(\ln n)^2\ln s+us(\ln n)^2)$.
Step 4d takes on the whole in the loop 4, $O(u^2(\ln n)^2\ln s)$ operations.
Hence the stated running time.
\end{proof}
\end{theorem}

Suppose now that $\mathcal{B}=\lbrace a,b\rbrace$ with $\p{ord}_n(a)$ and $\p{ord}_n(b)$ arbitrary.
Let $A=\p{ord}_n(a)$, $B=\p{ord}_n(b)$. For $s\in\mathbb{P}$, set $g_s=a^{\frac{A}{s^{v_s(A)}}}$ if 
$v_s(A)\ge v_s(B)$, else $g_s=b^{\frac{B}{s^{v_s(B)}}}$. We follow the procedures 
PH($n,a^{\frac{A}{s^{v_s(A)}}},b^{\frac{B}{s^{v_s(B)}}},s,v_s(A),v_s(B)$), $s$ running through the
set of primes dividing $(A,B)$. The group $\langle a,b\rangle_n$ is a direct sum of its $s$-primary parts 
$\langle a^{\frac{A}{s^{v_s(A)}}},b^{\frac{B}{s^{v_s(B)}}}\rangle_n$. Therefore, either a nontrivial factor of 
$n$ will be found, or $\langle a,b\rangle_n$ is cyclic, generated by $\prod\limits_{s\mid AB}g_s$.

\begin{corollary}
\label{cykliczna}
Assume we are given a subset $\mathcal{B}$ of $\mathbb{Z}_n^*$ and the complete factorization of all the integers $\p{ord}_n(b)$ for $b\in\mathcal{B}$. Then we can find a generator of 
$\langle\mathcal{B}\rangle_n$ or a nontrivial factor of $n$ in $O(\#\mathcal{B}\cdot(p+\ln n)(\ln n)^3)$ deterministic time, where $p$ is the greatest prime dividing the order of at least two distinct 
$b_1,b_2\in\mathcal{B}$ (put $p=0$ if there is no such prime).
\begin{proof}
Again, the correctness has been already discussed. We obtain the run-time bound by summing 
$(s+v_s(\varphi(n))\ln s)v_s(\varphi(n))(\ln n)^2$ over $s\mid\varphi(n)$, $s\le p$, and multiplying by 
$\#\mathcal{B}$.
\end{proof}
\end{corollary}

\begin{remark}
The number $p$ in the $O$ symbol above could be replaced by $\sqrt{p}\ln p$. To achieve this, one  uses
FFT techniques, well known from Pollard's \cite{pollard} or Strassen's \cite{strassen} algorithms for factoring
$n$ in $O(n^{\frac{1}{4}+\varepsilon})$ deterministic time. The hardest part of the PH() algorithm is finding
$j$, $0\le j<s$, such that $d_j>1$. The integer $j$ is of the form $j=j_0+j_1\lceil\sqrt{s}\rceil$ for some 
integers $j_0, j_1$, $0\le j_0, j_1<\lceil\sqrt{s}\rceil$. Let $a'=a^{s^{v-1}}$, $c'=c^{s^{v-m-2}}$.
We introduce the polynomial $h=\prod\limits_{0\le i_0<\lceil\sqrt{s}\rceil}(c'-a'^{i_0}X)$ and compute $h(a'^{i_1\lceil\sqrt{s}\rceil})$ for $i_1=0, 1, \ldots, \lceil\sqrt{s}\rceil-1$. By theorem 4 of \cite{turk}
it can be done in $O(\sqrt{s}(\ln s)^2(\ln n)^2)$ deterministic time. Next we find $j_1$ satisfying 
$(h(a'^{j_1\lceil\sqrt{s}\rceil}),n)>1$. Afterwards $j_0$ such that $(c'-a'^{j_0+j_1\lceil\sqrt{s}\rceil},n)>1$.
The computational cost of these last two steps is $O(\sqrt{s}(\ln n)^2)$, thus negligible.
\end{remark}

Turning back to our main question, we propose the following deterministic algorithm for
splitting $n$ given an integer $M$ as in the beginning of this section.
\newline
\newline
\noindent
Split($n,M,s_1,v_1,\ldots,s_t,v_t$)
\{$M=s_1^{v_1}\cdot\ldots\cdot s_t^{v_t}$ is the complete factorization of $M$\}
\begin{enumerate}
\item
For every $b\in\mathcal{B}$, compute $b^M$ modulo $n$, and:
\begin{enumerate}
\item
If $(b^M-1,n)=1$ then report failure and stop
\item
If $(b^M-1,n)<n$ then output this gcd and stop
\end{enumerate}
\item
Using the complete factorization of $M$, compute $\p{ord}_n(b)$ for each $b\in\mathcal{B}$
\item
For every $b\in\mathcal{B}$ and each prime $s\mid\p{ord}_n(b)$, compute
$(b^{\frac{\p{ord}_n(b)}{s}}-1,n)$. If one of these gcds is a nontrivial factor of $n$, then stop
\item
Using the algorithm associated with corollary \ref{cykliczna}, check whether 
$\langle\mathcal{B}\rangle_n$ is cyclic. If a nontrivial divisor of $n$ is found during these computations,
then stop
\item
State that $n$ is prime
\end{enumerate}

\begin{theorem}
\label{split}
Let ${\mathcal{B}}=\{2,3,\ldots,[(\ln n)^2]\}$,
$M=s_1^{v_1}\cdot\ldots\cdot s_t^{v_t}$ be the complete factorization of the integer $M$,
$s_0=\max\lbrace s\mid M:~\forall_{q\mid n}~s\mid q-1\rbrace\cup\{0\}$.
Suppose that $b^M\equiv 1 (n)$ for all $b\in\mathcal{B}$. Then the 
algorithm Split($n,M,s_1,v_1,\ldots,s_t,v_t$) finds a nontrivial divisor 
(or a proof of the primality) of $n$ in 
 $O((s_0\ln n+(\ln M)(\ln\ln M)+(\ln n)^2)(\ln n)^4)$ deterministic time.
\begin{proof}
For the correctness assume that we have reached step 5 of the algorithm. Step 3 implies that
$\mathcal{B}$ contains no Fermat-Euclid witness for $n$ and step 4 that 
$\langle\mathcal{B}\rangle_n$ is cyclic. Therefore $n$ is indeed prime in the light of the 
Fellows-Koblitz primality criterion. We proceed to the running time analysis. 
Step 1 requires $O((\ln M)(\ln n)^4)$ operations.
Step 2 can be done in $O((\ln M)(\ln\ln M)(\ln n)^4)$ time (see \cite{konyagin} - the analysis of the runtime of algorithm 3.1).
Step 3 costs $O${\Large$(\frac{(\ln n)^6}{\ln\ln n})$} operations. When we get to step 4, 
the exponent of $\langle\mathcal{B}\rangle_n$ divides $q-1$ for every prime factor $q$ of $n$. 
By corollary \ref{cykliczna}, the remaining computations thus take $O((s_0+\ln n)(\ln n)^5)$ time.
\end{proof}
\end{theorem}

There might be inputs $n$ for which the runtime of Split($n,M,s_1,v_1,\ldots,s_t,v_t$) is not
polynomial in $\ln n$ and $\ln M$, but it actually is if the integer $s_0$ defined in theorem \ref{split} is small, say bounded by a polynomial $B$ in $\ln n$. This is obviously satisfied whenever $n$ has a prime divisor $p$ such that $p-1$ is $B$-smooth.

\begin{corollary}[deterministic version of Pollard's $p-1$ algorithm]
\label{pollard}
Let $B\ge\ln n$.
\renewcommand{\labelenumi}{(\roman{enumi})}
\begin{enumerate}
\item
Assume $n$ has a prime divisor
$p$ such that $p-1$ is $B$-smooth. Then we can find a nontrivial divisor (or a proof of the primality)
of $n$ in $O(B(\ln n)^5)$ deterministic time.
\item
Suppose in addition that $n$ has at most one prime divisor $p$ such that $p-1$ is not $B$-smooth.
Then we can obtain the complete factorization of $n$, together with a primality proof for each of
the prime factors, in $O(B(\ln n)^6)$ deterministic time.
\end{enumerate}
\begin{proof}
Put $M=\prod\limits_{q\le B}q^{\bigl[\frac{\ln n}{\ln q}\bigr]}$ in theorem \ref{split}.
Part (i) follows, since 
$\ln M=O${\Large$(\frac{B}{\ln B}$}$\ln n${\Large$)$} and $\ln\ln M=O(\ln B)$.
For part (ii), simply consider the iteration of the
algorithm corresponding to part (i), combined with the Lenstra-Pomerance variant of the AKS
primality test \cite{isprime}, which runs in $O((\ln n)^6(\ln\ln n)^c)$ deterministic time for some 
constant $c$.
\end{proof}
\end{corollary}

Let us briefly compare the running times of the original Pollard $p-1$ algorithm with the new version.
The original algorithm finds a nontrivial divisor of $n$ in 
$O${\Large$(\frac{B}{\ln B}$}$(\ln n)^3${\Large$)$} random time under
the assumption of corollary \ref{pollard} (i).
Our deterministic version is slower (though not as much as we would expect)
and thus rather of theoretical than practical interest. 
\newline
Of course, the obtained running time bound of Split($n,M,s_1,v_1,\ldots,s_t,v_t$) is polynomial
in $\ln n$ and $\ln M$ for more inputs $n$ than those considered in corollary \ref{pollard},
with $B$ a polynomial in $\ln n$.
Let $D(n,u)=\max\limits_{q>(\ln n)^u}\#\{p\mid n:~q\mid p-1\}$, $u>0$.
We should expect that the integers $n$ for which $D(n,u)>1$ (with $u$ fixed) are rare. This is in fact true. 
We prove slightly more than needed to motivate the ideas of section \ref{dwa}.

\begin{theorem}
Let $l\in\mathbb{N}$. The number $B(x,u,l)$ of integers $n\le x$ such
that $D(n,u)>l$ is bounded above by 
$cx${\Large$\frac{2^{lu}(\ln\ln x)^{l+1}}{(\ln x)^{lu}}$}, where the constant $c$ does not depend 
upon $u$.
\begin{proof}
We have:
\begin{equation*}
B(x,u,l)\le\sqrt{x}+\sum\limits_{\sqrt{x}<n\le x}\sum\limits_{q>(\ln n)^u}
\sum\limits_{\substack{p_1<\ldots<p_{l+1}\\p_i\mid n\\p_i\equiv 1 (q)}} 1
\le\sqrt{x}+\sum\limits_{q>2^{-u}(\ln x)^u}\sum\limits_{\sqrt{x}<n\le x}
\sum\limits_{\substack{p_1<\ldots<p_{l+1}\\p_i\mid n\\p_i \equiv 1 (q)}} 1
\end{equation*}

\begin{align*}
\sum\limits_{n\le x}\sum\limits_{\substack{p_1<\ldots<p_{l+1}\\p_i\mid n\\p_i \equiv 1 (q)}} 1 &=
\sum\limits_{\substack{p_1<\ldots<p_{l+1}\le x\\p_i \equiv 1 (q)}}\biggl[\frac{x}{p_1\cdot\ldots\cdot p_{l+1}}\biggr]\le
x\sum\limits_{\substack{p_1<\ldots<p_{l+1}\le x\\p_i \equiv 1 (q)}}\frac{1}{p_1\cdot\ldots\cdot p_{l+1}}\\ &
\le x\biggl(\sum\limits_{\substack{p\le x\\p \equiv 1 (q)}}\frac{1}{p}\biggr)^{l+1}
\le\frac{c_1x(\ln\ln x)^{l+1}}{(q-1)^{l+1}},
\end{align*}

where the last inequality follows from the uniform bound

\begin{equation*}
\sum\limits_{\substack{p\le x\\p \equiv 1 (d)}}\frac{1}{p}\le\frac{c_0}{\varphi(d)}\ln\ln x
\end{equation*}

(use summation by parts and apply the Brun-Titchmarsh inequality). Hence

\begin{align*}
\sum\limits_{q>2^{-u}(\ln x)^u}\sum\limits_{\sqrt{x}<n\le x}
\sum\limits_{\substack{p_1<\ldots<p_{l+1}\\p_i\mid n\\p_i \equiv 1 (q)}} 1 &\le
c_1x(\ln\ln x)^{l+1}\sum\limits_{q>2^{-u}(\ln x)^u}\frac{1}{(q-1)^{l+1}}\\ &
\le c_2x\frac{2^{lu}(\ln\ln x)^{l+1}}{(\ln x)^{lu}}
\end{align*}

Thus
\begin{equation*}
B(x,u,l)\le c_3x\frac{2^{lu}(\ln\ln x)^{l+1}}{(\ln x)^{lu}}.
\end{equation*}

\end{proof}
\end{theorem}

\section{Generalization to the $p+1$ and other cyclotomic methods}
\label{cyklotomiczne}

Williams \cite{williams} designed a method of factoring analogous to Pollard's $p-1$ algorithm,
the $p+1$ method. It splits in random polynomial time integers $n$ having a prime divisor $p$
such that $p+1$ is smooth. Traditionally, it is described in terms of Lucas sequences, but the
analogy with the $p-1$ method becomes clear if one works, modulo $n$, in some quadratic
extension of $\mathbb{Z}$, as we will do. This section is mainly devoted to the proof of

\begin{theorem}
\label{williams}
Let $n$ and $m$ be odd, coprime integers, $n>2$, $m$ squarefree. Let $B\ge\ln n$.
Suppose that $n$ has a prime
factor $p$ such that $p+1$ is $B$-smooth and $\left(\frac{m}{p}\right)=-1$. Then we can find a nontrivial
divisor (or a proof of the primality) of $n$ in $O_{c,m}(B(\ln n)^{ch+3})$ deterministic time, 
where $h$ is the class number of $\mathbb{Q}(\sqrt{m})$ and $c$ is any constant greater than $4$.
\end{theorem}

The obtained derandomization of the $p+1$ algorithm is only partial, because of the requirement
$\left(\frac{m}{p}\right)=-1$, $m$ being fixed. We should therefore talk about deterministic $p+1$
methods (for varying $m$) instead of one deterministic $p+1$ algorithm. We need some auxiliary
results in the spirit of \cite{zralek}, the extension of the Pohlig-Hellman algorithm for the group $\mathbb{Z}_n[\sqrt{m}]^*$ to begin with.

\begin{theorem}
\label{PH2}
Suppose that $m$ mod $p$ is a quadratic nonresidue for some prime $p$ dividing $n$.
Let a subset $\mathcal{B}$ of $\mathbb{Z}_n[\sqrt{m}]^*$ and the complete factorization of all the integers 
$\p{ord}_{n,m}(b)$ for $b\in\mathcal{B}$ be given. Then a generator of $\langle\mathcal{B}\rangle_{n,m}$ or a nontrivial factor of $n$ can be computed in $O_m(\#\mathcal{B}\cdot(q+\ln n)(\ln n)^3)$ deterministic time, where $q$ is the greatest prime dividing the order of at least two distinct $b_1,b_2\in\mathcal{B}$ (set $q=0$ if there is no such prime).
\begin{proof}
As in corollary \ref{cykliczna} the argument reduces to the case of $\mathcal{B}=\{a,b\}$ with
$\p{ord}_{n,m}(a)$ and $\p{ord}_{n,m}(b)$ equal to the powers of some prime $s$, say 
$\p{ord}_{n,m}(a)=s^v$, $\p{ord}_{n,m}(b)=s^w$, $v\ge w$. Let $a^{s^{v-1}}=a_1+a_2\sqrt{m}$. We can 
also assume that $\p{ord}_{p,m}(a)=s^v$, for otherwise $(a_1-1,n)$ or $(a_2,n)$ would be a nontrivial divisor of $n$. The rest of the proof follows the lines of section \ref{jeden}, since $\mathbb{Z}_p[\sqrt{m}]^*$ is, by assumption, isomorphic to $\mathbb{F}_{p^2}^*$, hence cyclic.
\end{proof}
\end{theorem}

We introduce the standard integral basis of the ring of integers in $\mathbb{Q}(\sqrt{m})$, letting
$y=\sqrt{m}$ if $m\equiv 2,3 (4)$, and $y=\frac{1+\sqrt{m}}{2}$ if $m\equiv 1 (4)$. The next  theorem
is well known in the context of solving generalized Pell equations (norm equations in 
$\mathbb{Z}[y]$).

\begin{theorem}
There is an effective, positive constant $c_1$ depending upon $m$ and having the following property.
For any nonzero $a\in\mathbb{Z}[y]$, there exists $b\in\mathbb{Z}[y]$, $b=b_1+b_2y$, such that
$\frac{b}{a}\in\mathbb{Z}[y]^*$ and $|b_i|\le c_1\sqrt{|N(a)|}$, where $N(a)$ is the norm of $a$ and
$i=1,2$.
\end{theorem}

Finally, we formulate some kind of analogue of theorem \ref{konyaginpomerance} for the ring
$\mathbb{Z}[y]$.

\begin{theorem}
\label{baza}
Let $n$ be odd, $n>2$. Also, let $c>1$. Adopting the above notation, define 
\begin{equation*}
\mathcal{A}=\{a_1+a_2y:~|a_i|\le c_1(\ln n)^\frac{ch}{2}, 1\le i\le 2\},
\end{equation*}
\begin{equation*}
\mathcal{S}=\{v\cdot\alpha_1\cdot\ldots\cdot\alpha_t:~v\in\mathbb{Z}[y]^*, t\in\mathbb{N},
\alpha_i\in\mathcal{A}, 1\le i\le t\},
\end{equation*}
and $\pi_n:~\mathbb{Z}[y]\rightarrow\mathbb{Z}_n[\sqrt{m}]$ as
the obvious projection. Then $\#\pi_n(\mathcal{S})>n^{2-\frac{2}{c}-\varepsilon}+1$ for any $\varepsilon>0$
and $n\ge n_0$, $n_0=n_0(m, c, \varepsilon)$.
\begin{proof}
This is in fact a special case of lemma 3.5 from \cite{zralek}.
\end{proof}
\end{theorem}

Let $f_n$ be the endomorphism 
\begin{equation*}
a_1+a_2\sqrt{m}\mapsto(a_1-a_2\sqrt{m})(a_1+a_2\sqrt{m})^{-1}
\end{equation*}
of $\mathbb{Z}_n[\sqrt{m}]^*$. Let $\mathcal{U}$ be a set
of generators of the group of units $\mathbb{Z}[y]^*$, $\#\mathcal{U}\le 2$ ($\mathcal{U}$ could be
written explicitly), and let 
\begin{equation*}
\mathcal{B}_n=\pi_n(\mathcal{U}\cup\mathcal{A})\setminus\{0\}.
\end{equation*}
The algorithm below is a deterministic version of the $p+1$ factorization method. We justify the
correctness in the proof of theorem \ref{williams}. 
\\\\
Split2($n,c,m,M,s_1,v_1,\ldots,s_t,v_t$)
\{$c>4$, $M=s_1^{v_1}\cdot\ldots\cdot s_t^{v_t}$ is the complete factorization of $M$\}
\begin{enumerate}
\item
If $n$ is a nontrivial power $d^k$ then output $d$ and stop
\item
Let $n_0$ be as in theorem \ref{baza}, with $\varepsilon=\frac{1}{2}-\frac{2}{c}$.
If $n$ has a prime factor below $n_0$ then output such one and stop
\item
For each $a\in\mathcal{A}$, compute $N(a)$, let $\pi_n(a)=a_1+a_2\sqrt{m}$, and:
\begin{enumerate}
\item
If $1<(N(a),n)<n$ then output $(N(a),n)$ and stop
\item
If $n\mid N(a)$ then:
\begin{enumerate}
\item
If $(a_1,n)=1$ or $(a_2,n)=1$ then output failure and stop
\item
If $(a_1,n)<n$ then output this gcd and stop. Do the same with $(a_2,n)$
\end{enumerate}
\end{enumerate}
\item
For every $b\in f_n(\mathcal{B}_n)$, compute $b^M$, $b^M=b_1+b_2\sqrt{m}$, and:
\begin{enumerate}
\item
If $(b_1-1,n)=1$ or $(b_2,n)=1$ then report failure and stop
\item
If $(b_1-1,n)<n$ then output  this gcd and stop. Do the same with $(b_2,n)$
\end{enumerate}
\item
Using the complete factorization of $M$, compute $\p{ord}_{n,m}(b)$ for each
$b\in f_n(\mathcal{B}_n)$
\item
For every $b\in f_n(\mathcal{B}_n)$ and each prime
$s\mid\p{ord}_{n,m}(b)$, compute $b^{\frac{\p{ord}_{n,m}(b)}{s}}$,
$b^{\frac{\p{ord}_{n,m}(b)}{s}}=b_1+b_2\sqrt{m}$, and the gcds $(b_1-1,n)$, $(b_2,n)$.
If one of these gcds is a nontrivial factor of $n$ then stop
\item
Using the algorithm associated with theorem \ref{PH2}, check whether
$\langle f_n(\mathcal{B}_n)\rangle_{n,m}$ is cyclic. If a nontrivial divisor of $n$ is found during these
computations, then stop
\item
State that $n$ is prime
\end{enumerate}

\begin{proof}[Proof of theorem \ref{williams}]
Set $M=\prod\limits_{q\le B}q^{\bigl[\frac{\ln (n+1)}{\ln q}\bigr]}$.
First, we have to show that under our assumptions the algorithm will not report any failure. This could
happen only in step 3b(i) or 4a. Let $n\mid N(a)$ in step 3b. Then, in particular, $p\mid N(a)$ and thus
the element $\pi_p(a)$ is not invertible. Moreover, $\mathbb{Z}_p[\sqrt{m}]$ is isomorphic to the
field $\mathbb{F}_{p^2}$, since $\left(\frac{m}{p}\right)=-1$. 
We conclude that $\pi_p(a)$ must be zero, that is to say, $p\mid a_1$ and
$p\mid a_2$. Consequently, the algorithm cannot terminate in step 3b(i). Now let 
$b\in f_n(\mathcal{B}_n)$ in step 4. From step 3, $\mathcal{B}_n\subset\mathbb{Z}_n[\sqrt{m}]^*$,
so $b$ is correctly defined. The conjugation modulo $p$ is easily seen to be nothing but the
Frobenius map. The endomorphism $f_p$ thus raises the elements of $\mathbb{Z}_p[\sqrt{m}]^*$
to the power of $p-1$. As $M$ is a multiple of $p+1$, it follows that $b^M$ modulo $p$ must
be equal to $1$. Therefore no failure can be reported in step 4a.\\
Second, we should prove that $n$ is prime when step 8 is reached. Let us assume the contrary and
seek a contradiction. Denote by $q$ the least prime factor of $n$, and by $n'$ the squarefree part 
of $n$. Define $A$ as $\p{LCM}_{b\in f_{n'}(\mathcal{B}_{n'})}\p{ord}_{n',m}(b)$. 
From step 6, we have $A=\p{LCM}_{b\in f_q(\mathcal{B}_q)}\p{ord}_{q,m}(b)$.
By step 7, $\langle f_n(\mathcal{B}_n)\rangle_{n,m}$ is cyclic; so are its homomorphic images
$\langle f_{n'}(\mathcal{B}_{n'})\rangle_{n',m}$ and $\langle f_q(\mathcal{B}_q)\rangle_{q,m}$.
Thus 
\begin{equation*}
\#\langle f_{n'}(\mathcal{B}_{n'})\rangle_{n',m}=A=\#\langle f_q(\mathcal{B}_q)\rangle_{q,m}.
\end{equation*}
Hence 
\begin{equation*}
\#\langle f_{n'}(\mathcal{B}_{n'})\rangle_{n',m}\mid\# f_q(\mathbb{Z}_q[\sqrt{m}]^*).
\end{equation*}
Furthermore, $\#\langle f_{n'}(\mathcal{B}_{n'})\rangle_{n',m}\ge
\frac{\#\langle\mathcal{B}_{n'}\rangle_{n',m}}{\#\ker f_{n'}}$. From step 2, $n'\ge n_0$, which by
theorem \ref{baza} yields $\#\langle\mathcal{B}_{n'}\rangle_{n',m}>n'^{\frac{3}{2}}$. We will
evaluate $\#\ker f_{n'}$. Let $s$ be a prime dividing $n'$. If $\left(\frac{m}{s}\right)=-1$, then we
already know that $\#\ker f_s=s-1$. In the case when $\left(\frac{m}{s}\right)=1$, it is not hard to show that $f_s$ acts like the endomorphism $(a,b)\mapsto(ba^{-1},ab^{-1})$ of $\mathbb{Z}_s^*\oplus\mathbb{Z}_s^*$, and therefore $\#\ker f_s=s-1$. Consequently,
$\#\ker f_{n'}=\prod\limits_{s\mid n'}\#\ker f_s=\prod\limits_{s\mid n'}(s-1)$. Combining all the above,
we get 
\begin{equation*}
\#\langle f_{n'}(\mathcal{B}_{n'})\rangle_{n',m}>\# f_q(\mathbb{Z}_q[\sqrt{m}]^*)\cdot
\frac{q^\frac{3}{2}}{\# f_q(\mathbb{Z}_q[\sqrt{m}]^*)\cdot\#\ker f_q}\cdot(q^{-1}n')^\frac{1}{2}.
\end{equation*}
By the isomorphism theorem, 
$\# f_q(\mathbb{Z}_q[\sqrt{m}]^*)\cdot\#\ker f_q=\#\mathbb{Z}_q[\sqrt{m}]^*$, which is less than $q^2$.
From step 1, $q<n'^{\frac{1}{2}}$. Hence 
\begin{equation*}
\#\langle f_{n'}(\mathcal{B}_{n'})\rangle_{n',m}>\# f_q(\mathbb{Z}_q[\sqrt{m}]^*)\cdot 
q^{-1}n'^{\frac{1}{2}}>\# f_q(\mathbb{Z}_q[\sqrt{m}]^*). 
\end{equation*}
This contradicts the previously obtained inequality $\#\langle f_{n'}(\mathcal{B}_{n'})\rangle_{n',m}\le\# f_q(\mathbb{Z}_q[\sqrt{m}]^*)$.\\
The running time analysis is similar to that of algorithm Split; the role of the ``base set'' $\mathcal{B}$ is played here by $f_n(\mathcal{B}_n)$, whose cardinality is $O_m((\ln n)^{ch})$.
\end{proof}

Pollard's $p-1$ and William's $p+1$ algorithms are part of a family of factoring algorithms  called the cyclotomic methods. These were introduced by Bach and Shallit \cite{cyclo}, who proved, conditionally on the generalized Riemann hypothesis (GRH), the following. Let $\Phi_k$ be the $k$-th cyclotomic polynomial. An integer $n$ can be split in random polynomial time whenever $\Phi_k(p)$ is smooth for some prime $p$ dividing $n$, and integer $k$ polynomial in the size of $n$. If we fix $k$ and strengthen (reasonably, of course) the condition on $p$, it will eventually appear that neither GRH nor randomness
are necessary.

\begin{theorem}
\label{cyclotomic}
Let $F$ be a monic, irreducible polynomial of degree $k$ in $\mathbb{Z}[X]$, $k\ge 2$, such that the
extension $K$ of $\mathbb{Q}$, obtained by adjoining a root $\theta$ of $F$, is cyclic. Let $m\mid k$,
$m\ge 2$, and $B\ge\ln n$. Assume that $n$ is divisible by a prime $p$ with the property that 
$\Phi_m(p)$ is $B$-smooth and $F$ modulo $p$ is irreducible in $\mathbb{Z}_p[X]$. Then a nontrivial
factor (or a proof of the primality) of $n$ can be computed in $O_{c,\theta}(B(\ln n)^{ch+3})$ deterministic
time, where $h$ is the class number of $K$ and $c$ is any constant greater than $2k$. 
\end{theorem}

In the proof we will adopt two more pieces of notation. We will write $\mathcal{O}_K$ for the ring of integers of $K$. Furthermore, let $G$ be a group (written multiplicatively), 
$a\in G$, $\eta:~G\rightarrow G$, $V=\sum v_iX^i\in\mathbb{Z}[X]$. The expression $V(\eta)(a)$ will stand for $\prod\eta^i(a^{v_i})$, $\eta^i$ being the $i$-th iteration of $\eta$ ($\eta^0$ the identity).

\begin{proof}
There is no loss of generality in supposing that $n$ is coprime to the discriminant of $F$. 
The rings $\mathcal{O}_K/(n)$ and $\mathbb{Z}_n[\theta]$ are then isomorphic; we identify
them for convenience.
The Galois group of $K$ over $\mathbb{Q}$ consists of, say, $\psi_1, \ldots, \psi_k$. Denote by $\psi_{i,n}$ the automorphism of $\mathbb{Z}_n[\theta]$ induced by $\psi_i$. Let $f_{i,n}$ be the endomorphism
\begin{equation*} 
a\mapsto\prod\limits_{l\mid k,~l\ne m}\Phi_l(\psi_{i,n})(a)
\end{equation*}
of $\mathbb{Z}_n[\theta]^*$. The prime $p$ remains prime in $\mathcal{O}_K$;
let $\psi_j$ be the Frobenius over $(p)$. Then $f_{j,p}$ acts like 
$\mathbb{F}_{p^k}^*\ni a\mapsto a^\frac{p^k-1}{\Phi_m(p)}\in\mathbb{F}_{p^k}^*$.
Consequently, setting $M=\prod\limits_{q\le B}q^{\bigl[\frac{m\ln n}{\ln q}\bigr]}$ yields
$f_{j,p}(a)^M=1$ for any $a\in\mathbb{Z}_p[\theta]^*$. Up to now, we followed \cite{cyclo}. However, in order to compute deterministically a nontrivial factorization of $n$, we define
a ``base set'' of the form $f_{j,n}(\mathcal{B}_n)$. We do not know $j$ a priori, but in practice we can
work in turn with each endomorphism $f_{i,n}$, $i=1, \ldots, k$. An integral basis 
$\omega=(\omega_1, \ldots, \omega_k)$ of $\mathcal{O}_K$ and a finite set $\mathcal{U}$ of generators for $\mathcal{O}_K^*$ should be constructed independently of $n$, in a precomputation phase. Consider 
\begin{equation*}
\mathcal{A}=\{a_1\omega_1+\ldots+a_k\omega_k:~|a_i|\le c_1(\ln n)^\frac{ch}{k}, 1\le i\le k\},
\end{equation*}
where $c_1$ is the constant $c_3$ from theorem 3.4 of \cite{zralek}. Let $\pi_n$ be the projection
$\mathcal{O}_K\rightarrow\mathbb{Z}_n[\theta]$. Similarly to the proof of theorem \ref{williams},
we can assume that $\pi_n(\mathcal{U}\cup\mathcal{A})\setminus\{0\}\subset\mathbb{Z}_n[\theta]^*$
and put $\mathcal{B}_n=\pi_n(\mathcal{U}\cup\mathcal{A})\setminus\{0\}$. Again, let $q=p_{-}(n)$ and
let $n'$ be the squarefree part of $n$. Here also we can force $f_{j,n}(\mathcal{B}_n)^M=\{1\}$ and
further 
\begin{equation*}
\#\langle f_{j,n'}(\mathcal{B}_{n'})\rangle_{\mathbb{Z}_{n'}[\theta]^*}\mid\#f_{j,q}(\mathbb{Z}_q[\theta]^*).
\end{equation*}
This would follow from appropriate generalizations of steps 4-7
of algorithm Split2. Still, the extension of theorem \ref{baza} to $\mathcal{O}_K$ gives
\begin{equation*}
\#\langle f_{j,n'}(\mathcal{B}_{n'})\rangle_{\mathbb{Z}_{n'}[\theta]^*}\ge
\frac{\#\langle\mathcal{B}_{n'}\rangle_{\mathbb{Z}_{n'}[\theta]^*}}{\#\ker f_{j,n'}}>
\frac{n'^{k-\frac{k}{c}-\varepsilon}}{\prod\limits_{s\mid n'}\#\ker f_{j,s}}
\end{equation*}
if $\varepsilon>0$ and $n'$ exceeds
some constant $n_0$ independent of $n$. We have finally reached the interesting part of the proof,
which is bounding $\#\ker f_{j,s}$ for $s$ a prime factor of $n$. There are two cases to treat:\\
\renewcommand{\labelenumi}{(\roman{enumi})}
\begin{enumerate}
\item
$s$ stays prime in $\mathcal{O}_K$,
\item
$s$ splits in $\mathcal{O}_K$: $(s)=S_1\cdot\ldots\cdot S_e$, where the $S_i$ are distinct primes
of degree $d$, $d=\frac{k}{e}$, $e\ge 2$.
\end{enumerate}
Before we do this, note that $\psi_j$ has order $k$ (because $\psi_{j,p}$ has order $k$). Suppose that
(i) holds. The automorphism $\psi_j$ generates the Galois group of $K$ over $\mathbb{Q}$, isomorphic
by reduction modulo $s$ to the Galois group of $\mathcal{O}_K/(s)$ over $\mathbb{F}_s$, and so
$\psi_{j,s}$ is raising to the power of $s^r$ for some $r$ relatively prime to $k$, $r<k$.
Therefore $f_{j,s}$ acts as 
$\mathbb{F}_{s^k}^*\ni a\mapsto a^{\prod\limits_{l\mid k,~l\ne m}\Phi_l(s^r)}\in\mathbb{F}_{s^k}^*$.
It is easy to show that 
$\prod\limits_{l\mid k,~l\ne m}\Phi_l(X^r)=\prod\limits_{l\mid k,~l\ne m}\prod\limits_{t\mid r}\Phi_{tl}$.
This product is coprime to $\Phi_m$, since $m\mid k$. We apply B\'ezout's identity for polynomials
to see that $(\Phi_m(s),\prod\limits_{l\mid k,~l\ne m}\Phi_l(s^r))$ is bounded by a constant $c_2$
depending solely on $k$. Hence 
\begin{equation*}
\#\ker f_{j,s}=(s^k-1,\prod\limits_{l\mid k,~l\ne m}\Phi_l(s^r))\le
c_2\frac{s^k-1}{\Phi_m(s)}\le c_3s^{k-1},
\end{equation*}
where $c_3$ also depends only upon $k$.\\
Now assume that $s$ satisfies (ii). 
We want to bound the number of solutions $(a_1, \ldots, a_e)\in
(\mathcal{O}_K/S_1)^*\oplus\ldots\oplus(\mathcal{O}_K/S_e)^*$ to the equation 
$f_{j,s}(a_1, \ldots, a_e)=1$.
The automorphism $\psi_j$ acts on the set $\{S_1, \ldots, S_e\}$ as a cyclic permutation. 
In particular, $\psi_j^e$ generates the decomposition group of $S_1$, which is known to be isomorphic
(by reduction modulo $S_1$) to the Galois group of $\mathcal{O}_K/S_1$ over $\mathbb{F}_s$. Consequently, there is an $r$ coprime to $d$, such that $\psi_j^e(a)+S_1=a^{s^r}+S_1$ for every 
$a\in\mathcal{O}_K$. 
Thus $f_{j,s}(a_1, \ldots, a_e)+S_1$ is of the form 
$ba_1^{-1+\sum\limits_{0<i\le d-1}u_is^{ir}}$, with $b$ independent of $a_1$, and $u_i$ integers
depending just on $k$ and $m$. The $-1$ in the exponent of $a_1$ corresponds actually to the free
term of $\prod\limits_{l\mid k,~l\ne m}\Phi_l$ ($m\ge 2$). Since $(r,d)=1$, we have
$a_1^{-1+\sum\limits_{0<i\le d-1}u_is^{ir}}=a_1^{-1+\sum\limits_{0<i\le d-1}v_is^i}$, where the $v_i$
are a permutation of the $u_i$. In the field $\mathcal{O}_K/S_1$ there are at most 
$|-1+\sum\limits_{0<i\le d-1}v_is^i|$ solutions to the equation 
$ba_1^{-1+\sum\limits_{0<i\le d-1}v_is^i}=1$ of unknown $a_1$. Therefore
\begin{equation*}
\#\ker f_{j,s}\le(s^d-1)^{e-1}\cdot|-1+\sum\limits_{0<i\le d-1}v_is^i|\le c_4s^{k-1}
\end{equation*}
for a constant $c_4$ depending only upon $k$.\\
Proceeding along the same lines as the proof of theorem \ref{williams}, we get, if $\varepsilon>0$ and
$n'\ge n_0$,
the inequality 
\begin{equation*}
\#\langle f_{j,n'}(\mathcal{B}_{n'})\rangle_{\mathbb{Z}_{n'}[\theta]^*}>
\#f_{j,q}(\mathbb{Z}_q[\theta]^*)\cdot c_5q^{-1}n'^{1-\frac{k}{c}-\varepsilon},
\end{equation*}
where the (positive)
constant $c_5$ depends solely on $k$. Take $\varepsilon=\frac{1}{4}-\frac{k}{2c}$. Since
$\#\langle f_{j,n'}(\mathcal{B}_{n'})\rangle_{\mathbb{Z}_{n'}[\theta]^*}\le\#f_{j,q}(\mathbb{Z}_q[\theta]^*)$,
we conclude that $n$ is divisible by a prime less than $\max(n_0,c_5^\frac{4c}{2k-c})$, or $n$ is a
prime power.
\end{proof}

\begin{remark}
According to Frobenius' theorem, if $F$ is as in theorem \ref{cyclotomic}, then the set of primes $p$,
such that $F$ modulo $p$ is irreducible in $\mathbb{Z}_p[X]$, has density $\frac{\varphi(k)}{k}$.
This set consists in fact of primes lying in residue classes, which can be explicitly determined.
It suffices to express the root $\theta$ of $F$ as an element of a cyclotomic field (here we appeal to
the Kronecker-Weber theorem) and examine the order of the Frobenius automorphism in 
$\mathbb{Z}_p[\theta]$ (for $p$ not dividing the discriminant of $F$). As an example, $F=X^3-3X+1$
(a correct choice) is irreducible in $\mathbb{Z}_p[X]$ if and only if $p\equiv\pm 2~(9)$ or 
$p\equiv\pm 4~(9)$. We could thus reformulate theorem \ref{cyclotomic} in completely elementary terms
for specific polynomials $F$. We highly recommend that the reader interested in the theoretical setting of cyclotomic factoring algorithms, and willing to compare in detail our result with the classic method of Bach and Shallit, consult \cite{cyclo}. 
\end{remark}

\section{Some known reductions of factoring to computing $\varphi$}
\label{redukcje}

Taking $M=\varphi(n)$ in theorem \ref{probabilistycznypollard} we get the following classical
result.

\begin{theorem}[Rabin]
Given $\varphi(n)$ we can completely factor $n$ in $O((\ln n)^4)$ expected time.
\end{theorem}

For reasons already explained at the end of section \ref{standardowypollard}, substituting
$M=\varphi(n)$ also gives

\begin{theorem}[Miller]
If the ERH holds, then given $\varphi(n)$ we can completely factor $n$ in $O((\ln n)^6)$ deterministic time.
\end{theorem}

Define $G(n)$ as the least integer $m$ such that $\mathbb{Z}_n^*$ is generated by integers
less than or equal to $m$ and coprime to $n$. In \cite{burthe}, Burthe proved that 
$\frac{1}{x}\sum\limits_{n\le x}G(n)=O((\ln x)^{97})$. In particular, $G(n)<(\ln n)^{97+\varepsilon}$ for
almost all integers $n$. Now recall that any nonprincipal character modulo $n$ takes
a value different from $1$ for an integer less than or equal to $G(n)$. It follows by a similar argument
to the one used after theorem \ref{ankeny} that given $\varphi(n)$ we can completely factor $n$
in  $O((\ln n)^{101+\varepsilon})$ deterministic time for almost all $n$.  
\newline
While it is an open problem whether factoring unconditionally reduces in deterministic polynomial
time to computing  Euler's $\varphi$ function, for some integers such a reduction is particularly
easy. The simplest nontrivial examples are integers $n$ with exactly two prime factors.
Suppose first that $n=pq$. Then $p+q=n-\varphi(n)+1$. Given $\varphi(n)$ we compute the right-hand
side of this equality and find $p$ and $q$ by solving a quadratic equation. Now turn to the
general case $n=p^\alpha q^\beta$, say $p<q$. If $p\nmid q-1$, then 
$\frac{n}{(n,\varphi(n))}=pq$ and $\frac{\varphi(n)}{(n,\varphi(n))}=(p-1)(q-1)=\varphi(pq)$; thus the previous
method applies. If $p\mid q-1$, then $\frac{n}{(n,\varphi(n))}=q$ and therefore $q, \beta, \alpha, p$
will be obtained one after the other. 

Landau \cite{landau} showed that computing  the equal order factorization of any integer $n$,
that is, the sequence $n_i:=\prod\limits_{p:~v_p(n)=i}p$ ($i\ge 1$), can be done in deterministic polynomial
time given a ``$\varphi$-oracle'' (this oracle finds instantly the values of Euler's $\varphi$ function for 
$O(\ln n)$-bit inputs). In fact, if $\omega(n)\ge 3$, then $O(\Omega(n)(\ln n)^2)$ 
bit operations and at most $\omega(n)-2$ oracle calls (including $\varphi(n)$) are needed.
Notice that if $\omega(n_i)\le 2$ for all $i$, then the additional calls $\varphi(n_i)$ will lead to the complete factorization of $n$. For instance every integer $n=p^\alpha q^\beta s^\gamma$, where 
$p, q, s$ are
distinct primes and $\alpha, \beta, \gamma$ integers not all equal, can be, given $\varphi(n)$, completely factored in $O((\ln n)^3)$ deterministic time.

\section{Some subsets of the graph of $\varphi$ recognizable in deterministic polynomial time}
\label{dwa}

In section \ref{jeden} we have described in simple, arithmetic terms a set of integers of density $1$
in $\mathbb{N}$ (the set $\{n:~D(n,u)\le 1\}$ with $u$ fixed) whose elements $n$ are all factorable
in deterministic polynomial time if $\varphi(n)$ is given in a fully factored form. The ideas presented there are extended here to get a much more concrete result: exhibit a possibly large set of integers $n$
that are factorable in deterministic polynomial time given $\varphi(n)$ and only a part of its factorization,
which in turn can be obtained in polynomial time with the deterministic Pollard $p-1$ method.\\
Let $B$ and $\delta$ be positive real numbers. First define the following subsets of $\mathbb{P}$.
\begin{itemize}
\item
$\mathcal{P}_B$ is the set of primes $q$ such that $p-1$ is $B$-smooth for every prime $p$ dividing
$q-1$.
\item
$\mathcal{Q}_{B,\delta}$ is the set of primes $q$ such that the $B$-smooth part of $q-1$ is not less
than $q^\delta$.
\end{itemize}
Now consider, for $k$ an integer, $u, \delta, \eta$ positive real numbers, $\delta<1$, $\eta\le 1$,
the set $\mathcal{N}_{k,u,\delta,\eta}$ of integers that can be written in the form $n=n_1n_2n_3$,
where the $n_i$ are pairwise coprime, and: 
\begin{enumerate}
\item
$n_1$ has exactly $k$ distinct prime factors, all belonging to $\mathcal{P}_{(\ln n)^u}$.
\item
$n_2$ is a product of primes from $\mathcal{Q}_{(\ln n)^u,\delta}$.
\item
$n_3$ has at most two distinct prime factors. Furthermore, if $\omega(n_3)=2$ and $n_2\ne 1$,
then $p_{-}(n_2)>p_{-}(n_3)^\eta$.
\end{enumerate}
We will prove

\begin{theorem}
\label{graf phi}
Let $\mathcal{N}_{k,u,\delta,\eta}$ be as above. Given the pair $(n,\varphi(n))$, with 
$n\in\mathcal{N}_{k,u,\delta,\eta}$, we can completely factor $n$ in $O((\ln n)^C)$
deterministic time for some constant $C$ depending only on $k,u,\delta,\eta$. In particular, the set
$\{(n,\varphi(n)):~n\in\mathcal{N}_{k,u,\delta,\eta}\}$ is recognizable in deterministic polynomial time
($k,u,\delta,\eta$ being fixed).
\end{theorem}

We prepare the proof with some lemmas, keeping the notation of the theorem and assuming, without loss of generality, that $p_{-}(n)>(\ln n)^{\max(\frac{2}{\delta}, \frac{2+\eta}{\delta\eta}, k+3)}$.

\begin{lemma}
\label{lemme1}
Let $d$ be a factor of $n$, $M$ the $(\ln n)^u$-smooth part of $\varphi(n)$,
$\mathcal{B}=\{2, 3, \ldots, [(\ln d)^\frac{2}{\delta}]\}$ and $\mathcal{G}=\mathcal{B}^\frac{\varphi(n)}{M}$
modulo $d$. Assume that $d$ is divisible by two distinct primes $q_1, q_2$ from 
$\mathcal{Q}_{(\ln n)^u,\delta}$.
Then $\mathcal{G}$ contains a Fermat-Euclid witness for $d$ or 
$\langle\mathcal{G}\rangle_d$ is not cyclic.
\begin{proof}
Without loss of generality we let $q_1<q_2$.
Suppose, on the contrary, that there is no Fermat-Euclid witness for $d$ among the elements of 
$\mathcal{G}$
and that $\langle\mathcal{G}\rangle_d$ is cyclic. Then $\langle\mathcal{G}\rangle_{q_1q_2}$ is also
cyclic, as a homomorphic image of $\langle\mathcal{G}\rangle_d$, and so
$\#\langle\mathcal{G}\rangle_{q_1q_2}=\p{LCM}_{g\in\mathcal{G}}\p{ord}_{q_1q_2}(g)$. Moreover,
\begin{equation*}
\p{LCM}_{g\in\mathcal{G}}\p{ord}_{q_1q_2}(g)=\p{LCM}_{g\in\mathcal{G}}\p{ord}_d(g)=
\p{LCM}_{g\in\mathcal{G}}\p{ord}_{q_1}(g). 
\end{equation*}
Therefore $\#\langle\mathcal{G}\rangle_{q_1q_2}$
divides $(q_1-1,M)$, which equals, say $M_1$. We will show that 
$\#\langle\mathcal{G}\rangle_{q_1q_2}>M_1$ to derive a contradiction. Denote by $h$ the endomorphism raising every element of $\mathbb{Z}_{q_1q_2}^*$ to the power of 
$\frac{\varphi(n)}{M}$. We have 
$\langle\mathcal{G}\rangle_{q_1q_2}=h(\langle\mathcal{B}\rangle_{q_1q_2})$; hence
$\#\langle\mathcal{G}\rangle_{q_1q_2}\ge
\frac{\#\langle\mathcal{B}\rangle_{q_1q_2}}{\#\ker h}$. The numerator 
$\#\langle\mathcal{B}\rangle_{q_1q_2}\ge\psi(q_1q_2,(\ln q_1q_2)^\frac{2}{\delta})>
(q_1q_2)^{1-\frac{\delta}{2}}$. The denominator 
$\#\ker h=(q_1-1,\frac{\varphi(n)}{M})(q_2-1,\frac{\varphi(n)}{M})=\frac{(q_1-1)(q_2-1)}{M_1M_2}$,
where we let $M_2=(q_2-1,M)$. Also, $q_2\in\mathcal{Q}_{(\ln n)^u,\delta}$ and $q_2>q_1$; thus
$M_2\ge q_2^\delta>(q_1q_2)^\frac{\delta}{2}$. Putting all together gives
$\#\langle\mathcal{G}\rangle_{q_1q_2}>
M_1M_2\frac{(q_1q_2)^{1-\frac{\delta}{2}}}{(q_1-1)(q_2-1)}>M_1$, as required.
\end{proof}
\end{lemma}

\begin{lemma}
\label{lemme2}
Let $d$ be a factor of $n$, $M$ the $(\ln n)^u$-smooth part of $\varphi(n)$,
$\mathcal{B}'=\{2, 3, \ldots, [(\ln d)^\frac{2+\eta}{\delta\eta}]\}$ and
$\mathcal{G}'=\mathcal{B}{'}^\frac{\varphi(n)}{M}$ modulo $d$.
Suppose that $d$ is divisible by two distinct primes $p$ and $q$, 
$q\in\mathcal{Q}_{(\ln n)^u,\delta}$, $q>p^\eta$. Then $\mathcal{G}'$ contains a Fermat-Euclid witness for $d$ or $\langle\mathcal{G}'\rangle_d$ is not cyclic.
\begin{proof}
Suppose that neither element of $\mathcal{G}'$ is a Fermat-Euclid witness for $d$. We are to explain
why then  $\langle\mathcal{G}'\rangle_d$ cannot be cyclic. Let 
$A=\p{LCM}_{g\in\mathcal{G}'}\p{ord}_p(g)$. By assumption, $A$ also equals 
$\p{LCM}_{g\in\mathcal{G}'}\p{ord}_q(g)$, which is $\#\langle\mathcal{G}'\rangle_q$. Write $M_1$ for
the $(\ln n)^u$-smooth part of $q-1$. Similarly to the proof of lemma \ref{lemme1}, we obtain
\begin{equation*}
\#\langle\mathcal{G}'\rangle_q\ge M_1\frac{\#\langle\mathcal{B}'\rangle_q}{q-1}>
q^\frac{2\delta}{2+\eta}. 
\end{equation*}
Therefore $A> q^\frac{2\delta}{2+\eta}> p^\frac{2\delta\eta}{2+\eta}$.
Since $A$ divides $(p-1,M)$, it follows that $p\in\mathcal{Q}_{(\ln n)^u,\frac{2\delta\eta}{2+\eta}}$.
Furthermore, $q\in\mathcal{Q}_{(\ln n)^u,\delta}\subset\mathcal{Q}_{(\ln n)^u,\frac{2\delta\eta}{2+\eta}}$.
Replacing $\delta$ by $\frac{2\delta\eta}{2+\eta}$ in lemma \ref{lemme1}, we conclude that
$\langle\mathcal{G}'\rangle_d$ is not cyclic.
\end{proof}
\end{lemma}

\begin{lemma}
\label{lemme3}
Let $d$ be a factor of $n$, $M'=\prod p^{v_p(\varphi(n))}$, where the product ranges over the primes 
$p$ such that  $p-1$ is $(\ln n)^u$-smooth, $\mathcal{B}''=\{2, 3, \ldots, [(\ln d)^{k+3}]\}$. Assume that
$d$ has a prime divisor $q\in\mathcal{P}_{(\ln n)^u}$ and that $\omega(d)\le k+2$. Then one of the
following conditions holds.
\renewcommand{\labelenumi}{(\roman{enumi})}
\begin{enumerate}
\item
$1<(b^{M'}-1,d)<d$ for some $b\in\mathcal{B}''$.
\item
$b^{M'}\equiv 1 (d)$ for all $b\in\mathcal{B}''$ and $\mathcal{B}''$ contains a Fermat-Euclid
witness for $d$.
\item
$b^{M'}\equiv 1 (d)$ for every $b\in\mathcal{B}''$ and, setting 
$A=\p{LCM}_{b\in\mathcal{B}''}\p{ord}_d(b)$, we have $p_{-}(d)\equiv 1 (A)$, $A>d^\alpha$, with
$\alpha>\frac{1}{\omega(d)}-\frac{1}{\omega(d)^2}$.
\end{enumerate}   
\begin{proof}
The definitions of $M'$ and $q$ imply that $q-1\mid M'$. Consequently, $b^{M'}\equiv 1 (q)$ for
any $b\in\mathcal{B}''$. We shall therefore suppose that $b^{M'}\equiv 1 (d)$ for every
$b\in\mathcal{B}''$, that there is no Fermat-Euclid witness for $d$ in $\mathcal{B}''$, and verify the
properties of $A$. Under the latter assumption, $A\mid p-1$ for all primes $p$ dividing $d$,
for $p_{-}(d)$ in particular. That forces $(A,d)=1$ and so $(\#\langle\mathcal{B}''\rangle_d,d)=1$.
Hence $\langle\mathcal{B}''\rangle_d\le\bigoplus\limits_{p\mid d}C_{p-1}\le\mathbb{Z}_d^*$.
Therefore $\langle\mathcal{B}''\rangle_d$ contains, for each prime factor $q$ of $A$, at most
$\omega(d)$ linearly independent elements of order dividing $q^{v_q(A)}$. It follows that
$A^{\omega(d)}\ge \#\langle\mathcal{B}''\rangle_d$. Thus 
$A>\psi(d,(\ln d)^{k+3})^\frac{1}{\omega(d)}>d^\alpha$, where
$\alpha=\frac{1}{\omega(d)}(1-\frac{1}{k+3})$. Checking that 
$\alpha>\frac{1}{\omega(d)}-\frac{1}{\omega(d)^2}$ is straightforward.
\end{proof}
\end{lemma}

\begin{lemma}[Coppersmith et al.]
\label{lemme4}
Assume we are given integers $h>v>0$ and reals $\alpha, \beta, \gamma$ satisfying
$0<\alpha<1$, $0\le\beta<\gamma\le 1-\alpha$, $v(v+1)+\gamma h(h-1)-2(\alpha+\beta)vh<0$.
If $d$ is larger than some effectively computable constant, then all the divisors of $d$ of the form
$sx+r$, where  $0<r<s<d$, $s\ge d^\alpha$, $(r,s)=(s,d)=1$, $d^\beta\le x\le d^\gamma$,
can be found in deterministic polynomial time in $v, h, \ln d$. 
\end{lemma}

\begin{lemma}
\label{lemme5}
Let $r, s, d, l$ be integers and $\alpha$ a real number. Suppose that
$0<r<s<d$, $s\ge d^\alpha$, $(r,s)=(s,d)=1$, $\alpha>\frac{1}{l}-\frac{1}{l^2}$ and
$d$ is sufficiently large. Then all the divisors of $d$ of the form $sx+r$ and
less than $d^\frac{1}{l}$ can be found in $O_\varepsilon((\ln d)^3)$ deterministic time,
where $\varepsilon=\alpha-\frac{1}{l}+\frac{1}{l^2}$.
\begin{proof}
This is achieved by partitioning $[1,d^{\frac{1}{l}-\alpha}]$, the range of $x$, into intervals to which 
lemma \ref{lemme4} can be applied. We refer the reader to \cite{divisors} for the details of the algorithm.
For the running time, just follow closely the proof of lemma \ref{lemme4} therein.
\end{proof}
\end{lemma}

\begin{proof}[Proof of theorem \ref{graf phi}]
We describe an algorithm to compute the complete factorization of $n$.
\begin{enumerate}
\item
Let $L_1=\{n\}$, $L_2=\emptyset$
\item
Use the AKS primality test to check whether $L_1$ consists exclusively of prime numbers.
If so or $L_1=\emptyset$ then:
\begin{enumerate}
\item
If $L_2=\emptyset$ then output $n=\prod\limits_{p\in L_1}p^{v_p(n)}$ as the
complete factorization of $n$ and stop
\item
If $\#L_2>1$ then report failure and stop. In the contrary case, try to factor the
only element $m$ of $L_2$ into a product of two primes, $m=p^\alpha q^\beta$, assuming that
$\varphi(m)=\frac{\varphi(n)}{\prod\limits_{s\in L_1}s^{v_s(n)-1}(s-1)}$. If this works then output
$n=p^\alpha q^\beta \prod\limits_{s\in L_1}s^{v_s(n)}$ as the complete factorization of $n$ and stop.
Otherwise report failure and stop
\end{enumerate}
\item
Choose $d\in L_1\setminus\mathbb{P}$
\item
If $d$ is a prime power $p^\alpha$ then replace $d$ by $p$ in $L_1$. Return to step 2
\item
Attempt to split $d$ by means of the factoring algorithms corresponding evidently to lemmas
\ref{lemme1}, \ref{lemme2}, \ref{lemme3} and \ref{lemme5}. If this produces a nontrivial factorization
$d=d_1d_2$ then further apply a factor refinement procedure (cf. \cite{sigma}) to get a nontrivial
factorization $d=d_1'd_2'$ with $(d_1',d_2')=1$. Also, remove $d$ from $L_1$, adjoin $d_1',d_2'$
to $L_1$, and return to step 2
\item
Remove $d$ from $L_1$ and adjoin it to $L_2$. Return to step 2
\end{enumerate}
The algorithm obviously terminates. All we need to show is that when it does, $L_2=\emptyset$ or
$L_2=\{m\}$, with $\omega(m)=2$. Let $d$ be an integer chosen in step 3 of the algorithm, $d$ not
equal to a prime power. Then $d$ must have one of the following forms:
\renewcommand{\labelenumi}{(\roman{enumi})}
\begin{enumerate}
\item
$d$ divisible by two distinct primes from $\mathcal{Q}_{(\ln n)^u,\delta}$
\item
$d$ divisible by a prime from $\mathcal{P}_{(\ln n)^u}$, at most one prime from 
$\mathcal{Q}_{(\ln n)^u,\delta}$ and at most one prime factor of $n_3$
\item
$d$ divisible by a prime $q$ from $\mathcal{Q}_{(\ln n)^u,\delta}$ and the prime $p_{-}(n_3)$,
$\omega(n_3)=2$
\item
$d=p^{v_p(n)}q^{v_q(n)}$, where $q\in\mathcal{Q}_{(\ln n)^u,\delta}$, $p=p_{+}(n_3)$,
$\omega(n_3)=2$
\item
$d=n_3$, $\omega(n_3)=2$
\item
$d=n_3q^{v_q(n)}$, where $q\in\mathcal{Q}_{(\ln n)^u,\delta}$, $\omega(n_3)=1$
\end{enumerate}
The integer $d$ will be split in deterministic polynomial time:
\begin{itemize}
\item
In case (i) by lemma \ref{lemme1}.
\item
In case (ii) by lemmas \ref{lemme3} and \ref{lemme5}, since then $\omega(d)\le k+2$.
\item
In case (iii) by lemma \ref{lemme2}, because then $q>p_{-}(n_3)^\eta$.
\end{itemize}
Clearly, $d$ can be adjoined to $L_2$ only in cases (iv)-(vi), and if it is, no other element will.
 \end{proof}

\subsection*{Remarks}

In part 1 of the definition of $\mathcal{N}_{k,u,\delta,\eta}$, assuming that the prime factors of $n_1$
belong to $\mathcal{P}_{(\ln n)^u}$ is assuming that the part of $\varphi(n)$, which can be completely
factored in deterministic polynomial time with the $p-1$ method, is a multiple of
$\prod\limits_{q\mid n_1}(q-1)$. This assumption could be slightly relaxed by considering other
deterministic factoring methods, such as the $p+1$ methods of section \ref{cyklotomiczne}. Also, the condition $\omega(n_1)=k$ could be replaced by the weaker: if $q_1, \ldots, q_{k+1}$ are $k+1$
distinct primes dividing $n_1$, then the gcd of $q_1-1, \ldots, q_{k+1}-1$, is $(\ln n)^u$-smooth.\\\\
Primality testing is a special case of the problem of testing for membership in
$\{(n,\varphi(n)):~n\in\mathcal{N}_{k,u,\delta,\eta}\}$ or, more generally, in
$\{(n,\varphi(n)):~n\in\mathbb{N}\}$. Indeed, the set of primes can be identified with the subset
$\{(n,n-1):~n\in\mathbb{P}\}$ of the graph of $\varphi$. Before primality was known to be decidable
in deterministic polynomial time \cite{aks}, Konyagin and Pomerance \cite{konyagin} showed that
for any fixed, positive $u$ and $\delta$, the set $\{q:~q\in\mathcal{Q}_{(\ln q)^u,\delta}\}$ is recognizable
in deterministic polynomial time. Some of their ideas are used in this article, but in a more synthetic
way. 
\\\\
To conclude this section, we shall state without proof a result similar to theorem \ref{graf phi} for the sum of divisors function $\sigma$ (for a random polynomial time reduction of factoring to computing $\sigma$ 
cf. \cite{sigma}). Let $\mathcal{R}$ be a finite subset of $\mathbb{Z}$, and let $\mathcal{R}'$ be
the set of primes $q$ such that $\left(\frac{m}{q}\right)=-1$ for some $m\in\mathcal{R}$. Moreover, let
\begin{itemize}
\item
$\mathcal{P}_{\mathcal{R},B}$ be the subset of $\mathcal{R}'$ of such primes $q$ that for each prime
$p$ dividing $q+1$:
\begin{enumerate}
\item[]
$p-1$ is $B$-smooth or
\item[]
$p+1$ is $B$-smooth and $p\in\mathcal{R}'$
\end{enumerate}
\item
$\mathcal{Q}_{\mathcal{R},B,\delta}$ be the subset of $\mathcal{R}'$ of such primes $q$ that the
$B$-smooth part of $q+1$ is not less than $q^\delta$ 
\end{itemize}

To define $\mathcal{N}_{\mathcal{R}, k,u,\delta,\eta}$, replace in the definition of 
$\mathcal{N}_{k,u,\delta,\eta}$ the set $\mathcal{P}_{(\ln n)^u}$ by $\mathcal{P}_{\mathcal{R},(\ln n)^u}$,
the set $\mathcal{Q}_{(\ln n)^u,\delta}$ by $\mathcal{Q}_{\mathcal{R},(\ln n)^u,\delta}$, and add
a fourth condition:
\begin{enumerate}
\item[(4)]
$v_q(n_1n_2)$ is odd for all primes $q$ dividing $n_1n_2$
\end{enumerate}
Then the following analogue of theorem \ref{graf phi} holds.

\begin{theorem}
Given the pair $(n,\sigma(n))$, with $n\in\mathcal{N}_{\mathcal{R}, k,u,\delta,\eta}$, the complete
factorization of $n$ can be computed in $O((\ln n)^{C'})$ deterministic time,
where $C'$ is some constant depending only upon $\mathcal{R},k,u,\delta,\eta$. In particular, membership in 
\linebreak
$\{(n,\sigma(n)):~n\in\mathcal{N}_{\mathcal{R},k,u,\delta,\eta}\}$ is decidable in deterministic polynomial time (for $\mathcal{R},k,u,\delta,\eta$ fixed and $\mathcal{R}$ finite).
\end{theorem}

\section{A subexponential reduction of factoring to computing $\varphi$}

We shall abbreviate any expression of the form $\exp \Bigl((\ln x)^a (\ln\ln x)^{1-a}\Bigr)$ as $L(x,a)$. In this section we will first prove

\begin{theorem}
\label{faktoryzacja}
Suppose that $\varphi(n)$ is given in a completely factored form. Then the complete factorization 
of $n$ can be found in less than $L(n,\frac{1}{3})^{1+o(1)}$ deterministic time.
\end{theorem}

Then deduce

\begin{corollary}
\label{podwykladniczy}
Let $k=\min\{l\in\mathbb{N}:~\varphi^l(n)=1\}$. There is a deterministic algorithm that, given the sequence
$(n,\varphi(n),\varphi^2(n),\ldots,\varphi^k(n))$, outputs the complete factorization of $n$ in less than 
$L(n,\frac{1}{3})^{1+o(1)}$ time.
\begin{proof}
Let $1\le m\le k$.
Once we have found the complete factorization of $\varphi^m(n)$, we can compute, from theorem \ref{faktoryzacja}, the complete factorization of $\varphi^{m-1}(n)$ in less than $L(n,\frac{1}{3})^{1+o(1)}$ deterministic time. Since $\varphi^k(n)=1$ and $k\le 1+\log_2 n$, the corollary follows by induction. 
\end{proof}
\end{corollary}

In the proof of theorem \ref{faktoryzacja} we will exhibit a procedure that factors $n$ recursively,
that is, splits any previously computed, reducible divisor $d$ of $n$ further. Let $p=p_{-}(d)$. Additionally,
let $\alpha, \beta, \gamma$ be real numbers from the interval $(0,1)$, parameters to be optimally chosen.
Assume that $p>L(d,1-\alpha)$. Define $\mathcal{B}$ as $\{2,3,\ldots,[L(d,1-\alpha)]\}$, and
denote $\p{LCM}_{b\in\mathcal{B}}(\p{ord}_d(b))$ by $A$.

\begin{lemma}
\label{omega1}
Let $(1-\beta)(1-\gamma)\le 1-\alpha$. Suppose that $\mathcal{B}$ contains no Fermat-Euclid witness for 
$d$ and that $\omega(d)>${\Large $(\frac{\ln d}{\ln\ln d})^\beta$}. Then $p=mA+1$ for some
integer $m<L(d,(1-\beta)\gamma)$ if $p$ is sufficiently large.
\begin{proof}
We have
\begin{equation*}
L(p,1-\gamma)\le
\exp \Bigl((\frac{1}{\omega(d)}\ln d)^{1-\gamma}(\ln\ln d)^\gamma\Bigr)<
L(d,(1-\beta)(1-\gamma))\le L(d,1-\alpha),
\end{equation*}
where the last inequality holds if $d$ is large enough. Assume that $d$ is indeed such.
As $\mathcal{B}$ contains no Fermat-Euclid witness for $d$, it follows that
$A=\p{LCM}_{b\in\mathcal{B}}(\p{ord}_p(b))$. Consequently, 
$A=\#\langle\mathcal{B}\rangle_p\ge\psi(p,L(p,1-\gamma))$. By theorem \ref{canfield}, we obtain
$A\ge pL(p,\gamma)^{-1}$ if $p$ is sufficiently large.
We can write $p=mA+1$ for some $m\in\mathbb{N}$, because $A\mid p-1$.
Therefore $mA<p\le AL(p,\gamma)$.
Hence 
\begin{equation*}
m<L(p,\gamma)\le 
\exp \Bigl((\frac{1}{\omega(d)}\ln d)^{\gamma}(\ln\ln d)^{1-\gamma}\Bigr)<
L(d,(1-\beta)\gamma).
\end{equation*}
\end{proof}
\end{lemma}

\begin{lemma}
\label{omega2}
Let $\beta\le\frac{1}{2}$, $1-\beta\ge\alpha$. 
Assume that there is no Fermat-Euclid witness for $d$ in $\mathcal{B}$ and that
$\omega(d)\le${\Large $(\frac{\ln d}{\ln\ln d})^\beta$}.
Then $A^{\omega(d)+1}>d${\Large ${\omega(d)\choose [\omega(d)/2]}$} if $d$ is sufficiently large.
\begin{proof}
Just as in the proof of lemma \ref{lemme3}, we have 
$A^{\omega(d)}\ge\#\langle\mathcal{B}\rangle_d\ge\psi(d,L(d,1-\alpha))$. Hence 
$A^{\omega(d)+1}\ge\psi(d,L(d,1-\alpha))^{\frac{\omega(d)+1}{\omega(d)}}$. 
Let $-1<\varepsilon<-\alpha$.
It follows from theorem \ref{canfield} that 
$A^{\omega(d)+1}\ge 
d^{1+\frac{1}{\omega(d)}}L(d,\alpha)^{\varepsilon\frac{\omega(d)+1}{\omega(d)}}$ if $d$ is large enough.
It is sufficient to show that
$d^{\frac{1}{\omega(d)}}L(d,\alpha)^{\varepsilon\frac{\omega(d)+1}{\omega(d)}}>
{\omega(d)\choose [\omega(d)/2]}$ for large $d$.
This is clear when $\varepsilon\frac{\omega(d)+1}{\omega(d)}\le -1$, because then $\omega(d)$
is bounded from above. Suppose therefore that $\varepsilon\frac{\omega(d)+1}{\omega(d)}> -1$.
For sufficiently large $d$ we get  
\begin{align*}
d^{\frac{1}{\omega(d)}}L(d,\alpha)^{\varepsilon\frac{\omega(d)+1}{\omega(d)}}&\ge
L(d,1-\beta)L(d,\alpha)^{\varepsilon\frac{\omega(d)+1}{\omega(d)}}\ge
L(d,1-\beta)^{1+\varepsilon\frac{\omega(d)+1}{\omega(d)}}\\ &>
\exp \Bigl(\Bigl(\frac{\ln d}{\ln\ln d}\Bigr)^{\beta}\ln 2\Bigr)\ge
\exp (\omega(d)\ln 2)=2^{\omega(d)}\\ &>
{\omega(d)\choose [\omega(d)/2]}
\end{align*}
\end{proof}
\end{lemma}

The case $k=3$ of the ensuing lemma was proved in \cite{konyagin}.

\begin{lemma}
\label{poly}
Let $d=p_1^{e_1}\cdot\ldots\cdot p_k^{e_k}$.
Assume $A$ divides $p_i-1$ for $i=1,\ldots,k$; $p_i=b_iA+1$. Suppose in addition
that 
$A^{k+1}>{k\choose [k/2]}d$.
Write $d$ in base $A$: $d=1+a_1A+\ldots+a_kA^k$. Let $g=1+a_1X+\ldots+a_kX^k$. Then
$g=(b_1X+1)\cdot\ldots\cdot(b_kX+1)$ in $\mathbb{Z}[X]$. Furthermore, this factorization can
be obtained with the Hensel-Berlekamp algorithm in $O((\ln d)^5(\ln\ln d)^2)$ deterministic time.
\begin{proof}
We have $d=p_1^{e_1}\cdot\ldots\cdot p_k^{e_k}=(b_1A+1)^{e_1}\cdot\ldots\cdot (b_kA+1)^{e_k}$.
Since $A^{k+1}>d$, it follows that $e_1=\ldots=e_k=1$. Hence
\begin{equation*}
1+a_1A+\ldots+a_kA^k=(b_1A+1)\cdot\ldots\cdot (b_kA+1)
=1+\sum\limits_{j=1}^{k}\sigma_{k,j}(b_1,\ldots,b_k)A^j,
\end{equation*}
where $\sigma_{k,j}(b_1,\ldots,b_k)=\sum\limits_{1\le i_1<\ldots<i_j\le k}b_{i_1}\cdot\ldots\cdot b_{i_j}$.
It is therefore sufficient to show that $0\le\sigma_{k,j}(b_1,\ldots,b_k)<A$ for every $j$, $1\le j\le k$.
By assumption, $A^{k+1}>{k\choose [k/2]}d$ and thus
$b_1\cdot\ldots\cdot b_k {k\choose [k/2]}d<b_1\cdot\ldots\cdot b_k A^{k+1}<dA$. Hence
$b_1\cdot\ldots\cdot b_k {k\choose [k/2]}<A$ and it follows that
$0\le\sigma_{k,j}(b_1,\ldots,b_k)\le{k\choose j}b_1\cdot\ldots\cdot b_k
\le{k\choose [k/2]}b_1\cdot\ldots\cdot b_k<A$.
\newline
It remains to prove that $g$ can be completely factored in the stated time. We first need a ``small'' prime $p$ not dividing $a_k$ and such that $g_p$ is squarefree, $g_p$ being the reduction of $g$ modulo 
$p$. An upper bound for such a $p$ is given in \cite{lenstra} (3.9): $p=O(k\ln k+k\ln |g|)$, where
$|g|:=(1+\sum\limits_{i=1}^k a_i^2)^\frac{1}{2}$. 
Verifying that $p=O((\ln d)^2)$ is straightforward.
Let $\alpha=a_k^{-1} (p)$, $e=\lceil\frac{\ln d}{\ln p}\rceil$.
We factor completely $\alpha g_p$ with the Berlekamp algorithm in
$O(k(k+p)(k\ln p)^2)=O((\ln d)^5(\ln\ln d)^2)$ deterministic time (cf. theorem 7.4.5 of \cite{bachshallit}).
Then we lift this factorization to the factorization $\prod\limits_{1\le i\le k}(x+b_i^{-1})$ modulo $p^e$
with the Hensel algorithm in $O(ke(k\ln p)^2)=O((\ln d)^4(\ln\ln d)^2)$ deterministic time
(cf. theorem 7.7.2 of \cite{bachshallit}). Finally, we compute $b_i~(p^e)$ for every $i$. This finishes
the proof, as each $b_i$ is less than $p^e$. 
\end{proof}
\end{lemma}

\begin{proof}[Proof of theorem \ref{faktoryzacja}]
We find the complete factorization of $n$ using the algorithms associated with lemmas
\ref{omega1}, \ref{omega2} and \ref{poly}.
The running time bound of our recursive procedure is
obviously less than $L(n,\max(1-\alpha,(1-\beta)\gamma))^{1+o(1)}$. 
It remains to minimize $\max(1-\alpha,(1-\beta)\gamma)$ over the set
\begin{equation*}
\{ (\alpha,\beta,\gamma):~0<\alpha<1, 0<\beta\le\frac{1}{2}, 0<\gamma<1, 
1-\beta\ge\alpha, (1-\beta)(1-\gamma)\le 1-\alpha \}. 
\end{equation*}
Some easy calculations show
that the minimum is $\frac{1}{3}$, reached for $\alpha=\frac{2}{3}$, 
$\beta=\frac{1}{3}$, $\gamma=\frac{1}{2}$.
\end{proof}

\begin{remark}
The above method reduces the factorization of Carmichael numbers $n$ to the factorization
of $n-1$ in less than $L(n,\frac{1}{3})^{1+o(1)}$ deterministic time.
\end{remark}

\section*{Acknowledgements}

This paper contains part of the author's doctoral dissertation, written under the supervision of Dr. Jacek Pomyka\l a. It is a pleasure to thank him for all his help, encouragement and kindness.

\bibliographystyle{amsplain}

\end{document}